\documentstyle[leqno,11pt]{article}

\makeatletter
\@addtoreset{equation}{section}
\makeatother

\newcommand{\R}{I\hspace{-1.5mm}R}
\newcommand{\scriptR}{I\hspace{-0.9mm}R}

\newtheorem{Theorem}{Theorem}[section]

\newtheorem{Lemma}[Theorem]{Lemma}

\newtheorem{Remark}[Theorem]{Remark}
\date{}
\title{Uniqueness of solutions to 
Hamilton-Jacobi equations
arising in the Calculus of Variations}
\author{Gianni Dal Maso\thanks{SISSA, via Beirut 2-4, 
34014 Trieste, Italy} 
\hspace{.1cm}
 and H\'el\`ene Frankowska\thanks{CNRS, ERS2064, Centre de
Recherche Viabilit\'{e}, Jeux, Contr\^{o}le, Universit\'{e} de
Paris-Dauphine, 75775 Paris Cx 16, France} \hspace{.1cm} }

\begin{document}
\maketitle
\noindent
{\bf Abstract.}

\noindent
We prove the uniqueness of the viscosity solution to the Hamilton-Jacobi 
equation associated with a Bolza problem of the Calculus of Variations,
assuming that the Lagrangian is autonomous, continuous, superlinear, and
satisfies the usual convexity hypothesis. Under the same assumptions we
prove also the uniqueness, in a class of lower semicontinuous functions,
of a slightly different notion of solution, where classical derivatives are
replaced only by subdifferentials. These results follow from a new comparison
theorem for lower semicontinuous viscosity supersolutions of the
Hamilton-Jacobi equation, that is proved in the general case of lower
semicontinuous Lagrangians.

\vspace{ 5 mm}


\vspace{ 5 mm}

\noindent
{\bf AMS-MOS subject classification:} 49L20 (primary),
49L25 (secondary).

\section{Introduction}
Let us consider a {\em  Bolza problem\/} of the Calculus of
Variations
\begin{equation}\label{B}
\min \left\{ \int_{0}^{t} L (y (s),y' (s))ds \;+\;
\varphi (y(t)) : y \in W^{1,1} (0,t;\R^n),\; y(0)=x \right\},
\end{equation}
where the final cost $ \varphi \colon  \R^n \mapsto \R_+ \cup  \{ + \infty
\}$ is {\em  lower semicontinuous\/} and the Lagrangian
$L\colon  \R^n
\times \R^n \mapsto \R_+$  is {\em locally bounded\/} and {\em lower
semicontinuous\/}. We assume also that $L (x, \cdot )$
is {\em  convex\/} for every  $x \in \R^n$ and that the following
Tonelli type
{\em coercivity assumption\/} is satisfied: there exists
a function $\Theta \colon  \R^n \mapsto \R_+$ such that
\begin{equation} \label{T}
 \lim_{|u| \rightarrow \infty } \frac{ \Theta (u)}{|u|} =+
\infty \,, \qquad\forall \; (x,u)\in \R^n\times
\R^n,\;\; L (x,u) \geq \Theta ( u )\,.
\end{equation}
These assumptions guarantee the existence of absolutely continuous
minimizers (see, e.g., \cite{ces83}). The classical Lagrange problem (with the fixed final condition
$y (t)=z$) may be reduced to the form  (\ref{B}) by simply setting $
\varphi
(z):=0$ and $ \varphi := + \infty $ elsewhere.

The {\em value function\/}
$V \colon  \R_+ \times  \R^n \mapsto \R_+ \cup \{+ \infty \}$ for the Bolza
problem (\ref{B}) is
defined by
\begin{displaymath}
 V (t,x):= \min \left\{ \int_{0}^{t} L (y (s),y' (s))ds \;+\; \varphi
(y (t)) : y\in W^{1,1} (0,t;\R^n),\; y(0)=x \right\}.
\end{displaymath}
Under our assumptions on  $L$ and  $ \varphi $ the value function is  {\em
lower
semicontinuous\/} on  $\R_+ \times \R^n$ and 
{\em
locally Lipschitz\/} on
$\R_+^{ \star } \times \R^n$, where $\R_+^{ \star }:=
\R_+ \setminus \{0\}$. Moreover it satisfies the  {\em initial condition\/}
\begin{displaymath}
 \forall \; x \in \R^n,\;\; \liminf_{h \rightarrow 0+,\; y \rightarrow x} V
(h,y)= \varphi  (x) \,.
\end{displaymath}

If  $V$ is smooth and  $L$ is continuous, then  $V$ is a classical solution
to the {\em  Hamilton-Jacobi equation\/}
\begin{equation}\label{H-J}
\cases{V_t+H(x,-V_x)=0& in $\R_+^{\star} \times \R^n$,
\cr\cr
V (0, \cdot )= \varphi& in $\R^n$,
\cr}
\end{equation}
with Hamiltonian $H\colon  \R^n
\times \R^n \mapsto \R$ defined by
\begin{equation}\label{H}
 H (x,p) := \sup_{u \in \scriptR^n} \left( \left\langle   p,u\right\rangle
 -L(x,u) \right),
\end{equation}
where  $ \left\langle   \cdot , \cdot \right\rangle $ denotes the scalar
product in  $\R^n$. 
In other words, $H (x, \cdot )$ is the Legendre-Fenchel conjugate of $L (x,
\cdot )$.

It is  well known  that  (\ref{H-J}) may not have smooth solutions,
even if  $H$ and $\varphi $ are smooth. To overcome this difficulty,
different notions of
generalized
solutions have been proposed. The notion of viscosity
solution can be introduced by means of subdifferentials and
superdifferentials. 

We recall that the
{\em
subdifferential\/}  $ \partial _-W (x)$ of 
  a function $W\colon \R^n \mapsto \R \cup \{ +
\infty \}$ at a point  $x \in \mbox{\rm dom}(W)$ is defined by
\[
\partial_- W(x):=\left\{
p\in \R^n\;:\; \liminf_{y\to x}\frac{ W(y)-
W(x)-\langle p,y-x\rangle}{|y-x|} \geq 0\right\} .
\]
while the {\em superdifferential \/}  $ \partial _+ W (x)$ is defined by
$ \partial _+W (x):=- \partial _- (-W) (x)$. 

A  {\em continuous function\/} 
$W\colon  \R_+^{ \star } \times \R^n \mapsto \R$ is said to be a  {\em viscosity
solution to the Hamilton-Jacobi equation  (\ref{H-J})\/} if the following
conditions are satisfied:
\begin{equation}
 \label{hj1}
\forall \;  (t,x)\in\R_+^{ \star } \times \R^n, \; \forall \; (p_t,p_x) \in
\partial _-W (t,x),\;\; p_t + H (x,-p_x) \geq 0\,,
\end{equation}
\begin{equation}
 \label{hj2a}
\forall \; (t,x) \in\R_+^{ \star } \times \R^n, \; \forall \; (p_t,p_x) \in
\partial _+W (t,x),\;\; p_t + H (x,-p_x) \leq 0\,.
\end{equation}
In \cite{cl2} and \cite{cl2a} the uniqueness of a
bounded uniformly continuous viscosity solution of (\ref{H-J}) is proved
under some assumptions on  $H$, which imply in particular that
$H$ is continuous and  $H ( \cdot ,p)$ is uniformly continuous for
every $p \in\R^n$ . In 
\cite[Theorem 4.5]{dm-fra} we proved the following result, that can be
applied also to unbounded solutions.
\begin{Theorem}
\label{dm1.1}
Assume that $L$ is continuous. Let  $W\colon\R_+^{ \star } \times \R^n
\mapsto \R_+$ be a locally Lipschitz viscosity solution of  (\ref{H-J})
which satisfies the initial condition
\begin{equation}\label{hj2}
\forall \; x \in \R^n,\;\; \liminf_{h \rightarrow 0+,y \rightarrow x} W
(h,y)= \varphi  (x) \,.
\end{equation}
Then  $W=V$ on  $\R_+^{ \star } \times \R^n $.
\end{Theorem} 

To describe the new uniqueness results, we introduce the notion of {\em
contingent directional derivatives\/} of a function $W\colon \R^n \mapsto \R
\cup \{+ \infty \}$. These are defined, for every   $x \in \mbox{\rm
dom}(W)$ and for every  $u \in \R^n$, by
\begin{equation}\label{dirder}
D _{\uparrow } W (x) (u) :=
\liminf_{ \textstyle {h \rightarrow 0+ \atop v \rightarrow
u}} \frac{W (x+hv) - W (x)}{h} \;\;.
\end{equation}

The main result of this paper is the following theorem, which shows that
the value function  $V$ is the unique viscosity solution in the larger
class of continuous functions whose contingent derivatives satisfy the
following very weak assumptions:
\begin{equation}
 \label{new}
\forall \;  (t,x) \in \mbox{\rm dom}(W),\;\; t >0,\;\;D _{\uparrow }W (t,x)
(0,0) =0 \,,
\end{equation}
\begin{equation}
 \label{new1}
\forall \;  (t,x) \in \mbox{\rm dom}(W),\;\; t >0,\;\;\exists \;
u\in\R^{n}, \;\;D _{\uparrow }W (t,x) (-1,u) <+\infty\,.
\end{equation}
\begin{Theorem} \label{visc}
Assume that $L$ is continuous.
 Let $W \colon  \R_+^{ \star }
\times \R^n \mapsto \R_+ $ be a continuous viscosity solution of 
(\ref{H-J}) which satisfies  (\ref{new}),  (\ref{new1}), and the initial
condition (\ref{hj2}). Then $W=V$ on  $\R_+^{ \star }
\times \R^n$.
\end{Theorem}

In  \cite[Theorem 4.4]{dm-fra} we considered also a different notion of
generalized solution and we proved the following uniqueness result in the
class of locally Lipschitz functions.
\begin{Theorem}
Assume that $L$ is continuous.
 Let $W \colon  \R_+^{ \star }
\times \R^n \mapsto \R_+ $ be a locally Lipschitz function which satisfies
the initial condition  (\ref{hj2}) and solves the Hamilton-Jacobi equation 
(\ref{H-J}) in the following sense :
\begin{equation}
\label{dm-fra3}
 \forall \; (t,x) \in\R_+^{ \star } \times \R^n, \; \forall \; (p_t,p_x)
\in \partial _-W (t,x),\;\; p_t + H (x,-p_x) = 0 \,.
\end{equation}
 Then $W=V$ on $\R_+^{ \star } \times \R^n$.
\end{Theorem} 

In this paper we shall prove the following result, which provides
uniqueness in the larger class of lower semicontinuous functions $W$
satisfying  (\ref{new}) and  (\ref{new1}).
\begin{Theorem} \label{dm-fra4}
Assume that $L$ is continuous.
 Let $W \colon  \R_+^{ \star }
\times \R^n \mapsto \R_+ \cup \{+ \infty \}$ be a lower semicontinuous function which
satisfies the initial condition  (\ref{hj2}), the technical conditions 
(\ref{new}) and  (\ref{new1}), and
solves the Hamilton-Jacobi equation  (\ref{H-J}) in the sense of 
(\ref{dm-fra3}).
Then $W=V$ on $\R_+^{ \star } \times \R^n$.
\end{Theorem} 

In the proof of Theorem \ref{dm-fra4} we use the following comparison
result, which follows immediately from  \cite[Theorem 4.1, Remark 4.2 and
Proposition  4.3]{dm-fra}.
\begin{Theorem} \label{dm-fra5}
Assume that $L$ is continuous.
 Let $W \colon  \R_+^{ \star }
\times \R^n \mapsto \R_+ \cup \{ +
\infty \}$ be a lower semicontinuous function which
satisfies the initial condition  (\ref{hj2}). Suppose that  $W$ is a
subsolution of the Hamilton-Jacobi equation   (\ref{H-J}) in the following
sense:
\begin{equation} \label{viscsubsol}
\forall \; (t,x) \in\R_+^{ \star } \times \R^n, \; \forall \; (p_t,p_x)
\in \partial _-W (t,x),\;\; p_t + H (x,-p_x) \leq  0
\end{equation}
Then $W \leq V$ on $\R_+^{ \star } \times \R^n$.
\end{Theorem}

In the proof of Theorems \ref{visc} and \ref{dm-fra5} we use also a very
general comparison result (Theorem \ref{dm-fra1.6})
 for lower semicontinuous viscosity
supersolutions of the Hamilton-Jacobi equation  (\ref{H-J}). To our
knowledge the strongest result in this direction, dealing with possibly
discontinuous Lagrangians, is the following theorem proved in 
\cite[Theorem
5.1]{dm-fra}, where the notion of supersolution is given by using
contingent inequalities.
\begin{Theorem}
\label{dm-fra1.5}
 Let $W \colon  \R_+
\times \R^n \mapsto \R_+ \cup \{ +
\infty \} $ be a lower semicontinuous function which
satisfies the initial condition  $W (0, \cdot )= \varphi $. Suppose that 
$W$ is a supersolution of (\ref{H-J}) in the following sense :
$$ \forall \; (t,x) \in \mbox{\rm dom}(W),\;t>0, \; \exists \; u \in
\R^n,\;\; D _{\uparrow }W (t,x) (-1,u) \leq  -L (x,u) \,. $$
Then $W \geq V$ on $\R_+ \times \R^n$.
\end{Theorem} 

The comparison result for viscosity supersolutions we are going to prove
is the following theorem, where we need the additional assumptions 
(\ref{new}) and  (\ref{new1}).

\begin{Theorem}
\label{dm-fra1.6}
 Let $W \colon  \R_+
\times \R^n \mapsto \R_+ \cup \{ +
\infty \}$ be a lower semicontinuous function which
satisfies (\ref{new}) and  (\ref{new1}).
Suppose that $W (0, \cdot )= \varphi $ and  that  $W$ is a viscosity
supersolution of the Hamilton-Jacobi equation (\ref{H-J}), i.e., $W$
satisfies  (\ref{hj1}).
Then $W \geq V$ on $\R_+ \times \R^n$.
\end{Theorem} 

\section{Preliminaries}

Let $K\subset \R^n$ be a nonempty subset and $x\in K$. The
contingent cone $T_K(x)$ to $K$ at $x$ is defined by
\[
v\in T_K(x) \; \Longleftrightarrow  \;\liminf_{h\to 0+}
\frac{\mbox{\rm dist}(x+hv,K)}{h}=0 \, .
\]
The negative polar cone $T^-$ to a subset $T\subset \R^n$ is given by
\[
T^- :=\{v\in \R^n\,:\;\forall \; w\in T\,,\;\langle v,w \rangle\leq 0\} \, .
\]
When $K$ is convex , then $\left[T_K (x) \right]^-$
coincides with the the usual normal cone $N_{K} (x)$ of convex analysis.

The epigraph ${\cal E}pi(W)$ of a function
$W\colon \R^n \mapsto \R \cup \{ +
\infty \}$ is defined by
\begin{displaymath}
 {\cal E}pi(W) := \left\{ (x,r)\in \R^n\times\R : r \geq W (x) \right\} .
\end{displaymath}
We shall
need the following
version of Rockafellar's result (see \cite{Rockafellar}).
\begin{Lemma}\label{Rockafellar}
Let $x\in \mbox{\rm dom}(W)$ and let
$(p,0)\in\left[T_{{\cal E}pi(W)}(x,W(x))\right]^-$ be such that
$p\neq 0$. Then there exist
$x_\varepsilon $ converging to $x$ (as $ \varepsilon  \rightarrow 0+$) and
\begin{displaymath}
 (p_\varepsilon,q_\varepsilon)\in\left[T
_{{\cal E}pi(W)}(x_\varepsilon,W(x_\varepsilon ))\right]^-
\end{displaymath}
converging to $ (p,0)$ as  $ \varepsilon  \rightarrow 0+$
such that $q_{ \varepsilon }<0$ for every $\varepsilon>0$.
\end{Lemma}

A closed subset $K$ of $\R^n$ is called a viability domain
of a set-valued map $G\colon K \leadsto  \R^n$ if for every $x\in K$
\[
G(x)\cap T_K(x)\neq\emptyset\,.
\]
The following statement summarizes several versions of the viability
theorem (see~\cite{Aubin}).

\begin{Theorem}[Viability]\label{viability}
Let $K\subset \R^n$ be a closed set and let $G\colon K\leadsto \R^n$ be
 an upper semicontinuous
set-valued map with compact convex values. 
The following conditions are equivalent:
\begin{description}
\item[\mbox{\rm (a)}] $K$ is a viability domain of $G$;
\item[\mbox{\rm (b)}]  $G(x)\cap\overline {co}\,T_K(x)\neq\emptyset$
for every $x\in K$;
\item[\mbox{\rm (c)}]  for every $x\in K$ there exist
$ \varepsilon >0$ and a solution $y\colon {[0,\varepsilon[}\mapsto K$
to the
Cauchy problem
\end{description}
\begin{equation}\label{diG}
\cases{
y'(t)\in G(y(t)) \,,
\cr
y(0)=x \,.&
\cr}
\end{equation}
\end{Theorem}

The equivalence (a) $\Longleftrightarrow$ (b)
was proved in \cite{Ushakov}. This proof was simplified in
\cite[page~85]{Aubin}. The fact that (a)
$\Longleftrightarrow$ (c) was first proved by
Bebernes and Schuur in \cite{beb}. A proof can be found
in \cite{AubinCellina} or \cite{Aubin}.

\vspace{ 5 mm}
The next theorem allows to deal with some unbounded
set-valued maps with closed convex values. As usual,  $B_R$ denotes
the closed ball with centre $0$ and radius $R$.

\begin{Theorem}\label{viability0}
Let $K \subset \R^n$ be a closed set and let
$G\colon K\leadsto \R^n$ be an upper semicontinuous set-valued map
with closed convex values. We assume that for every $x \in K$ there
exists $R>0$ such that
for all small $h>0$
\begin{displaymath}
 \mbox{\rm dist}  (x+h G(x),K)= \mbox{\rm dist} (x+h (G (x) \cap 
 B_{R}),K) \,.
\end{displaymath}
Then the following statements are equivalent:
\begin{description}
\item[\mbox{\rm (a)}]  $K$ is a viability domain of $G$;
\item[\mbox{\rm (b)}] for every $x\in K$ and
for every
$p \in \left[   T_K (x)\right]^-$
we have $\displaystyle
\inf_{u \in G(x)}\left\langle  p,u \right\rangle \leq 0$.
\end{description}
\end{Theorem}

Theorem~\ref{viability0} is a direct consequence of the following more 
technical result, which will be crucial in the proof of 
Theorem~\ref{dm-fra1.6}.

\begin{Theorem}\label{viability1}
Let $K \subset \R^n$ be a closed 
set, let $x\in K$, and let $G\colon K\leadsto \R^n$ be a set-valued map 
with non-empty closed convex values. 
Assume that there exists $R>0$ such that
for all small $h>0$
\begin{equation}\label{dist}
 \mbox{\rm dist}  (x+h G(x),K)= \mbox{\rm dist} (x+h (G (x) \cap B_R),K) 
 \,.
\end{equation}
Assume also that the support function, defined by
$$
 \sigma (y,p) := \sup_{u \in G (y)} 
 \left\langle   p,u\right\rangle \,,
$$
satisfies the following upper semicontinuity condition at $x$:
for every $\varepsilon>0$ there exists a 
neighbourhood $U$ of $x$ such that
\begin{equation}\label{usc}
\sigma(x,p)+\varepsilon|p|> \sigma(y,p)
\end{equation}
for every $y\in U\cap K$ and for every $p$ in the set
$$
N_{G(x)}(B_{R}):=\{p\in\R^{n}:
\exists\; u\in G(x)\cap B_{R},\  p\in N_{G(x)}(u)\}\,.
$$
Finally, assume that for every $y\in K$ in a suitable neighbourhood of 
$x$ we have
\begin{equation}\label{geq0}
\sigma (y,-p) \geq 0 
\end{equation}
for every
$p \in \left[   T_K (y)\right]^-$ such that 
$ -p \in N_{G (x)} (B_{R}) $.
Then $G(x)\cap T_K(x) \cap B_R\neq\emptyset$.
\end{Theorem}

{}From the proof given below it follows that the same result
holds if $\left[   T_K (y)\right]^-$ is replaced by the set of
proximal normals to $K$ at~$y$.

\vspace{ 5 mm}

{\bf Proof} --- \hspace{ 2 mm} Let us
define the  function $g \colon  \R_{+}\mapsto \R_+$ by
\begin{displaymath}
 g (h) := \frac{1}{2}\mbox{\rm dist} (x+hG (x),K)^2.
\end{displaymath}
Observe that  $g$ is locally Lipschitz around zero and $g (0)=0$.
For all small $h >0$ let us consider
$u_{h} \in G (x) \cap B _R$ and $x_h \in K$ such that $\mbox{\rm dist}
(x+hG
(x),K)=
|x+hu_h -x_h|$.
Then $x_h\rightarrow x$ when $h \rightarrow 0+$. 

Consider $h>0$ such that
$g' (h)$ does exist and fix any $u \in G (x)$.
Since
$G(x)$ is convex, for all nonnegative $h, \;k$ we have $ (h+k) G (x)=h G
(x)
+k G (x)$. 
 Therefore
\begin{displaymath}
 g' (h) \leq \frac{1}{2}\lim_{k \rightarrow 0+} \frac{|x+ hu_{h}+ku 
-x_h|^2
- |x+ hu_h -x_h|^2 }{k} = \left\langle p_h, u  
\right\rangle ,
\end{displaymath}
where $p_h :=x+ hu_h -x_h$.
Consequently
\begin{equation}\label{gprime}
 g' (h) \leq \inf_{u \in G (x)}\left\langle   p_h, u\right\rangle =
-\sup_{u \in G (x)} \left\langle   -p_h,u\right\rangle  = -\sigma  
(x,-p_h)\, .
\end{equation}

As $x_h$ is a point of $K$ with minimum distance from $x+hu_h$, the
vector $p_h$ is a proximal normal to $K$, therefore it
belongs to $\left[   T_K (x_h)\right]^-$.
On the other hand, $u_h$ is the
point of $G (x)$ with minimum distance from  $(x_{h}-x)/h$. Thus 
$-p_h \in N_{G (x)} (u_h)\subset N_{G (x)} (B_{R})$. By (\ref{geq0}) 
we have
\begin{equation}\label{geq00}
-\sigma(x_{h},-p_{h})\le 0\,.
\end{equation}

By the uniform upper semicontinuity (\ref{usc}) of $\sigma$
for every $\varepsilon>0$ there exists $h_{\varepsilon}>0$
such that for $0<h<h_{\varepsilon}$
\begin{equation}\label{gpr}
-\sigma  \left( x, -p_h\right) \leq -\sigma \left( 
x_h,-p_h\right) + |p_h| \varepsilon\leq
-\sigma \left( x_h,-p_h\right) + hR\varepsilon\, ,
\end{equation}
where the last inequality follows from the 
fact that $|p_h| \leq hR$.
{}From (\ref{gprime}), (\ref{geq00}), and (\ref{gpr}) we obtain that
$
 g' (h) \leq hR \varepsilon
$
for every $h\in(0,h_{\varepsilon})$ at which the derivative $ g' (h)$ 
exists.

Integrating $g'$ we deduce that for $0<h<h_{\varepsilon}$
\begin{displaymath}
 g (h) \leq  \frac{Rh^2}{2}\varepsilon\,.
\end{displaymath}
This implies that for $0<h<h_{\varepsilon}$
\begin{displaymath}
\frac{\mbox{\rm dist} (x+hu_{h},K) }h=
\frac{\mbox{\rm dist} (x+hG (x),K)}h  \leq \sqrt{R\varepsilon}\, .
\end{displaymath}
Let $h_i \rightarrow 0+$ be a sequence such that  $u_{h_i}$ converges to
some $u \in G (x)$. From the very definition of contingent cone we deduce
that $u \in T_K(x)$.
$\; \; \Box$

\section{Proof of the comparison theorem}
This section is devoted to the proof of the new comparison theorem
for viscosity supersolutions.

\vspace{ 5 mm}

{\bf Proof of Theorem \ref{dm-fra1.6}} --- \hspace{ 2 mm} 
We first claim that $H$ is locally bounded. Indeed for all $x\in \R^{n}$ 
we have $H (x,p) \geq -L (x,0)$, thus $H$ is locally bounded from below. 
On the other hand,
$H (x,p) \leq \Theta^{*} (p)$, where  $\Theta^{*}$ denotes the
Legendre-Fenchel conjugate of  $\Theta$. Since the function  $\Theta$ has
a superlinear 
growth,  the convex function $\Theta^{ * }$ takes 
only finite values, so it is locally bounded. This shows that 
$H$ is also locally bounded from above. 

Consequently, 
the function $H(x,\cdot)$ is 
locally Lipschitz with respect to $p$, locally uniformly with respect to
$x$.
By (\ref{T}) $H(\cdot,p)$ is upper semicontinuous with respect to $x$.
These two properties together imply that for every $x\in \R^{n}$, for every
$M>0$, and for every $\varepsilon>0$ there exists a neighbourhood $U$ 
of $x$ such that
\begin{equation}\label{uscH}
\forall\; y\in U\,,\ \forall\; p\in B_{M}\,, \quad
H(x,p)+\varepsilon > H(y,p) \,.
\end{equation}

Let us define the set-valued map $G\colon  \R^n
\leadsto  \R \times \R^n \times \R$ with closed convex values by
\begin{displaymath}
  G (x):= \left\{  (-1,u,-L (x,u) - \rho ) :  \rho  \geq 0,\; u \in
\R^n \right\} .
\end{displaymath}
Let $K:= {\cal E}pi(W)$, let $ (t,x) \in \mbox{\rm dom}(W)$  with 
$t>0$,
and let $z:=(t,x,W(t,x))$. We want 
to show that all assumptions of Theorem \ref{viability1}
are satisfied (here
$z$ plays the role of $x$). To prove (\ref{dist})
we show that there exists $R>0$ 
such that for all small $h>0$
\begin{eqnarray*}
 &\exists \; u \in \R^n,\;\; | (-1,u,-L (x,u))| \leq R\, ,&
\\
 &\mbox{\rm dist}\left( (t,x,W(t,x)) + h (-1,u,-L (x,u)),K \right)=
 \mbox{\rm dist}\left(  (t,x,W(t,x))+h G (x),K \right)\,. &
\end{eqnarray*}
As $W\ge0$, it is easy to see that for every $h>0$  
there exist $ u_h \in \R^n$ and 
$z_h=(t_h,x_h, r_h)\in 
\R \times \R^n \times \R$, 
with $r_h\ge W(t_h,x_h)$,
such that
\begin{eqnarray*}
&|(t,x,W(t,x)) + h (-1,u_h, -L (x,u_h) ) - (t_h,x_h, r_h)| &
\\
&= \mbox{\rm dist} ((t,x,W(t,x))+hG (x), K)
\leq h (L (x,0)+1)\,. &
\end{eqnarray*}
If $W(t,x)-hL (x,u_h)-r_h>0$, by increasing $r_h$ we could make 
$W(t,x)-hL (x,u_h)-r_h=0$, which would contradict the definition of
distance.
Therefore $W(t,x)-hL (x,u_h)-r_h\le 0$, hence
\begin{displaymath}
hL (x,u_h) +W (t_h,x_h) -W (t,x) \le
hL (x,u_h) +r_h -W (t,x) \leq h(L (x,0)+1)\,.
\end{displaymath}

We claim that $u_h$ is bounded for all $h>0$ small enough.
Assume by contradiction that there exist $h_i  \rightarrow  0+$ such that
$|u_{h_i}| \rightarrow \infty $.
\vspace{ 3 mm}

{\it Case 1.\/} Assume first that for a subsequence, still 
 denoted by  $h_{i}$, we have $h_i|u_{h_i}| \geq c$ for some $c>0$.
Since $W(t_h,x_h)  \ge0$, we have
\begin{displaymath}
 h_iL (x,u_{h_i}) -W (t,x) \leq h_i (L (x,0) +1)\,;
\end{displaymath}
dividing by $h_i|u_{h_i}|$ and taking the limit we get
\begin{displaymath}
 \limsup_{i \rightarrow \infty } \frac{L (x,u_{h_i})}{|u_{h_i}|} <
 +\infty \, ,
\end{displaymath}
which contradicts  (\ref{T}).

\vspace{ 3 mm}

{\it Case 2.\/} It remains to consider the case $h_i |u_{h_i}|
\rightarrow 0$. Since
$$ \max \left\{ |t_h +h -t | ,\;\;  |x_h - x - hu_h| \right\} \leq h(L
(x,0)+1),$$
 we deduce that for some $v \in \R^n$ and for a subsequence, still 
 denoted by  $h_i$,
 we have
\begin{displaymath}
\lim_{i\rightarrow  \infty } \frac{t_{h_i} - t}{h_i |u_{h_i}|} =0
\qquad\mbox{and}\qquad \lim_{i \rightarrow \infty } \frac{x_{h_i} -x}{h_i
|u_{h_i}|}=
v\,.
\end{displaymath}
Furthermore,
\begin{displaymath}
h_iL (x,u_{h_i}) +W (t_{h_i},x_{h_i})  -W (t,x) \leq h_i(L (x,0)+1)\,;
\end{displaymath}
dividing by $h_i|u_{h_i}|$ and taking the limit yields
\begin{displaymath}
 \lim_{i\rightarrow \infty } \frac{W (t_{h_i},x_{h_i})  -W (t,x)}{h_i
|u_{h_i}|} = - \infty \,.
\end{displaymath}
Hence $D _{\uparrow }W (t,x)  (0,v)=- \infty $, which implies
$D _{\uparrow }W (t,x)  (0,0)=- \infty $ (see \cite{af90sva}).
This contradicts (\ref{new}) and completes the proof of our
claim.

\vspace{ 3 mm}

Thus $|u_h| $ is uniformly bounded when $h>0$ is small.
As $L$ is locally bounded, there exists $R>0$ such that
$| (-1,u_h, -L (x,u_h))| \leq R$ for all small $h>0$.

Observe next that, if  $(p_t,p_x,q) \in N_{G (x)} (B _R)$, then $q\ge 
0$. Moreover, if $q=0$, then $p_{x}=0$; if $q>0$,
then there exists $u\in  B _R$ such that
$ p_x/q  \in \partial_u L (x,u)$, the subdifferential of 
$L(x,\cdot)$ at $u$. As $L$ is locally bounded, this implies that
there exists a 
constant $M$ such that $p_x/q\in B_M$
for every
$(p_t,p_x,q) \in N_{G (x)} (B _R)$ with 
$q> 0$. 

To prove (\ref{usc}) it is enough to show that
for every $\varepsilon>0$ there exists a neighbourhood $U$ of $x$ 
such that
\begin{equation}\label{usc1}
\sup_{u\in\scriptR^{n}} (-p_t + \left\langle   p_x,u\right\rangle -qL (x,u)) +
\varepsilon |(p_t,p_x,q)|
\end{equation}
$$
> \sup_{u\in\scriptR^{n}} (-p_t + \left\langle   p_x,u\right\rangle -qL (y,u))
$$
for every $y\in U$ and for every $(p_t,p_x,q)\in N_{G (x)} (B _R)$.
If $q=0$, then $p_{x}=0$ and (\ref{usc1}) is trivial. If $q>0$, then
(\ref{usc1})
can be written as
$$
-p_t + qH(x,\frac{p_x}{q}) +
\varepsilon |(p_t,p_x,q)|> 
-p_t + qH(y,\frac{p_x}{q})\,,
$$
which follows easily from (\ref{uscH}).

Let us check that (\ref{geq0}) holds true. Fix
$(s,y,r) \in {\cal E}pi (W)$, with $s>0$, and
$ (p_t,p_x,q) \in \left[  T_{ {\cal E}pi (W)} (s,y,r)\right]^- $. 
Since $(0,0,1)\in T_{ {\cal E}pi (W)} (s,y,r)$, we have $ q \leq 0$.
Therefore (\ref{geq0}) is equivalent to
\begin{equation}\label{geq01}
\sup_{u\in\scriptR^{n}} (p_t + \left\langle 
 - p_x,u\right\rangle +qL (y,u))\ge 0 \,.
\end{equation}
If $q<0$, then $(p_t/|q|,p_x/|q|,-1)\in \left[  T_{ {\cal E}pi (W)} 
(s,y,r)\right]^-$, hence 
$ (p_t/|q|,p_x/|q|) \in \partial_-W (s,y)$ (see \cite[page~249]{af90sva}) 
and we deduce (\ref{geq01}) from  (\ref{hj1}). 
If $q=0$ and $p_x\neq0$, then the supremum in (\ref{geq01}) is 
$+\infty$. If $q=0$ and $p_x=0$, then
$(p_t,0,0)\in \left[  T_{ {\cal E}pi (W)} 
(s,y,r)\right]^-$. By (\ref{new1}) there exists $u\in \R^{n}$ such 
that $z:=D _{\uparrow }W (s,y) (-1,u) <+\infty$.
As ${\cal E}pi (D _{\uparrow } W (s,y) ( \cdot,\cdot ) ) =
 T_{{\cal E}pi (W)} (s,y,W (s,y)) $, the vector 
 $(-1,u, z)$ belongs to
 $T_{{\cal E}pi (W)} (s,y,W (s,y))$, which is contained in 
 $T_{ {\cal E}pi (W)}(s,y,r)$. By the definition of $\left[  T_{ {\cal E}pi
(W)} 
(s,y,r)\right]^-$ we obtain $-p_t\le 0$, which yields (\ref{geq01})
when $q=0$ and $p_x=0$.

{}From Theorem \ref{viability1} we deduce that
\begin{displaymath}
 G (x) \cap \left(  T_{ {\cal E}pi (W)} (t,x,W (t,x)) \right) \ne 
 \emptyset\,.
\end{displaymath}
As ${\cal E}pi (D _{\uparrow } W (t,x) ( \cdot,\cdot ) ) =
 T_{{\cal E}pi (W)} (t,x,W (t,x))$, we obtain that
\begin{displaymath}
 \exists \; u \in \R^n,\;  D _{\uparrow }W (t,x) (-1,u) \leq -L (x,u) .
\end{displaymath}
This and Theorem  \ref{dm-fra1.5} imply that $W \geq V$ on  $ \R_+ \times
\R^n$. $\; \; \Box$

\begin{Remark}\label{rem}
{ \rm {}From the proof of Theorem \ref{dm-fra1.6} we see that condition
(\ref{new1}) yields, for $t>0$,
$$
 (p_t,0,0) \in \left[  T_{ {\cal E}pi (W)} (t,x,W(t,x))\right]^-
 \; \Longrightarrow  \;
 p_t  \ge 0  \,.
$$
If, in addition, the subsolution inequality
(\ref{viscsubsol}) is satisfied, and 
$$(p_t,0,0) \in \left[  T_{ {\cal E}pi (W)} (t,x,W(t,x))\right]^-,$$
then by Rockafellar's Lemma \ref{Rockafellar} there exist
$ (t_i,x_i) \rightarrow (t,x)$
and $(p^i_{t},p^i_{x},q_i) \in
\left[   T_{{\cal E}pi (W)} (t_i,x_i, W(t_i,x_i)) \right]^-$,
with $q_i <0$, such that $(p^i_{t},p^i_{x},q_i) \rightarrow (p_t,0,0)$.
Then $(-p^i_t/q_i,-p^i_x/q_i) \in \partial_- W (t_i,x_i)$ and from
(\ref{viscsubsol}) we obtain $-p^i_t/q_i - L(x_i,0)\le 0$, hence
$p^i_t + q_i L(x_i,0)\le 0$. Taking the limit we get 
$p_t \le 0$. This shows that (\ref{new1}) and (\ref{viscsubsol}) together imply 
the following geometric condition for $t>0$:
$$
 (p_t,0,0) \in \left[  T_{ {\cal E}pi (W)} (t,x,W(t,x))\right]^-
 \; \Longrightarrow  \;
 p_t  = 0  \,.
$$}
\end{Remark}

\section{Proofs of the uniqueness results}
We begin with the proof of the uniqueness theorem for viscsity solutions.

\vspace{ 5 mm}

{\bf Proof of Theorem \ref{visc}} --- \hspace{ 2 mm} By Theorem 
\ref{dm-fra1.6}
we have
 $W \geq V$. Recall
that the hypograph of $W$ is defined by
\begin{displaymath}
 {\cal H}yp (W):=\{ (t,x,r) \in \R_+^{\star} \times \R^n \times
\R : r \leq W (t,x)\}\, .
\end{displaymath}
Define the closed set
$$
K := {\cal H}yp (W) \cup ( \R_- \times \R^n \times \R )\, .
$$
We claim that for all $(t,x,r) \in K$
\begin{equation} \label{cd12}
\forall \; u \in \R^n,\;\; (-1,u, -L (x,u) ) \in \overline {co}
\;T_{K} (t,x,r)\, .
\end{equation}
It is enough to prove it in the case $t>0$ and $r=W (t,x)$, since the other
cases are evident.
Fix $u \in \R^n$. Then by (\ref{hj2a})
\begin{equation} \label{cd13}
 \forall \; (p_t,p_x) \in \partial_+ W (t,x),\; p_t + \left\langle
-p_x,u\right\rangle -L (x,u)  \leq 0 \,.
\end{equation}

We want to prove that
\begin{equation} \label{cd14}
\forall\; (-p_t,-p_x,q) \in \left[   T_{ {\cal H}yp (W)} (t,x,W
(t,x))\right]^- \;\;
p_t + \left\langle   -p_x,u\right\rangle -qL (x,u)  \leq 0\, .
\end{equation}
When $q >0$ we have
$(p_t/q,p_x/q) \in \partial_+ W (t,x)$, thus (\ref{cd14})
follows from (\ref{cd13}).
By Lemma  \ref{Rockafellar}, applied to $-W$, if
$(0,0,0)\ne (-p_t,-p_x,0) \in
\left[ T_{{\cal H}yp (W)}
(t,x, W (t,x)) \right] ^-$, then for some $ (t_i,x_i) \rightarrow (t,x)$
and $(-p^i_{t},-p^i_{x},q_i) \in
\left[   T_{{\cal H}yp (W)} (t_i,x_i, W(t_i,x_i)) \right]^-$,
with $q_i >0$, we have $(p^i_{t},p^i_{x},q_i) \rightarrow (p_t,p_x,0)$.
So
\begin{displaymath}
 p^i_{t} + \langle   -p^i_{x},u\rangle -q_iL (x_i,u)  \leq 0 .
\end{displaymath}
Taking the limit we get
$ p_t + \left\langle   -p_x,u\right\rangle  \leq 0$, which concludes
the proof of (\ref{cd14}).

By the separation theorem, (\ref{cd12}) follows from
(\ref{cd14}).  Since the lower set-valued  limit of contingent cones is
equal
to Clarke's tangent cone  (see for instance  \cite{af90sva}),
from the continuity of $L$ we deduce that, for all
$(t,x) \in \R_+^{ \star } \times \R^n,$
\begin{displaymath}
 (-1,u,-L (x,u)) \in C_{ {\cal H}yp(W)}  (t,x,W (t,x)) .
\end{displaymath}
Fix $\varepsilon >0$.
Then it is not difficult to check that
\begin{displaymath}
\forall \; u \in \R^n,\;\;  (-1,u,-L (x,u) -\varepsilon )
\in \mbox{\rm Int} \left(
C_{ {\cal H}yp(W)}
 (t,x,W (t,x)) \right) .
\end{displaymath}
By \cite[Proposition 13, p. 425]{AE} this yields
\begin{equation} \label{B++}
 \liminf_{h \rightarrow 0+, \,v \rightarrow u} \frac{W (t -h, x+hv)-W
(t,x)}{h}  \geq -L (x,u) -\varepsilon  .
\end{equation}

We  have to show that for all $ t>0, \; x \in \R^n$, $V (t,x) \geq
W (t,x)$. Let $y$  be a minimizer of the Bolza
problem (\ref{B}). Since $L$ is
continuous, $y' \in L^{ \infty } (0,t;\R^{n})$ by \cite{ambr}.
Consider a sequence of continuous functions $u_i\colon
[ 0,t] \mapsto \R^n$ which is bounded in $L^{ \infty } (0,t;\R^{n})$
and converges to $y'$ almost
everywhere in $[0,t]$, and let $t_i \rightarrow 0+$ and $x_i \rightarrow y
(t)$
be such that $W (t_i,x_i)\rightarrow \varphi (y (t))$.
Define
\begin{displaymath}
 y_i (s) :=  x_i - \int_{s}^{t-t_i} u_i ( \tau )d \tau,\;\; \forall \; s \in 
[0,t-t_i] .
\end{displaymath}
and $y_i (s):=x_i$ for $s > t-t_i$.
Then $y_i$ converges to $y$ uniformly in $[0,t]$.
Fix  $i$ and set $ \psi (s):=W (t-s,y_i (s))$ for $s \in {[ 0,t-t_i]}$. By
(\ref{B++})
for every $s\in {[ 0,t-t_i[}$ we have
\begin{equation}
 \label{coine+}
\limsup_{h\to 0+ }  \frac{\psi  (s+h) -\psi  (s)}h \geq -L (y_i (s),  u_i (s))-\varepsilon .
\end{equation}
Consider  the system
\begin{displaymath}
\cases{
(\alpha ' (s),   z' (s) ) = (1, -L (  y_i (s),  u_i (s))
-\varepsilon), \;\; s  \geq 0\, &
\cr\cr
(\alpha  (0), z (0))  = (0, W (t, y_i (0)))\, .&
\cr}
\end{displaymath}
where we have set $u_i (s):=u_i (t)$ for all  $s \geq t-t_i$.
It has the unique solution
$$( \alpha  (s),z (s)) :=  \left( s, W (t, y _i (0)) -
\int_{0}^{s} L (y_i (
\tau ), u_i ( \tau )) d \tau - \varepsilon s \right) .$$
According to Theorem \ref{viability} and  (\ref{coine+}),
this solution is
viable in ${\cal H}yp ( \psi ) \cup ({[t-t_i, +\infty [} \times \R)$, i.e.,
for
all $s \in {[0,t-t_i]}$ we have
 $ ( \alpha  (s),z (s)) \in {\cal H}yp ( \psi )$. Thus for all $s \in
{[0,t-t_i]}$
$$W (t-s,y_i (s)) \geq  W (t,y_i (0)) - \int_{0}^{s} L ( y_i (
\tau ),   u_i ( \tau ))d \tau  - \varepsilon s\,.$$
In particular
\begin{displaymath}
 W (t, y_i (0)) \leq W (t_i,x_i) + \int_{0}^{t-t_i}  L ( y_i (
\tau ),   u_i ( \tau ))d \tau  + \varepsilon (t-t_i)\,.
\end{displaymath}
Since the functions $L$ and $W$ are continuous, 
taking the limit  we
get
\begin{displaymath}
 W (t,x)- \varepsilon  t \leq \varphi (y (t)) + \int_{0}^{t} L ( y ( \tau
),   y' ( \tau ))d \tau V (t,x) \,.
\end{displaymath}
Finally, as $ \varepsilon  \rightarrow 0+$ we obtain $W (t,x) \leq V
(t,x)$.
$\; \; \Box$

\vspace{ 5 mm}

{\bf Proof of Theorem \ref{dm-fra4}} --- \hspace{ 2 mm} 
The inequality $W\ge V$ is proved in Theorem \ref{dm-fra1.6},
while the inequality $W\le V$ follows from Theorem \ref{dm-fra5}.
$\; \; \Box$

\vspace{ 5 mm}

\end{document}